\documentclass[12pt,leqno]{amsart}
\usepackage{latexsym,amssymb,amsmath,multicol, rotating,lscape}
\usepackage{letltxmacro}

\usepackage{listings}

\usepackage{tikz}
\usetikzlibrary{calc,decorations.markings}

 2

\newtheorem{thm}{Theorem}[section]
\newtheorem{lem}[thm]{Lemma}
\newtheorem{prop}[thm]{Proposition}
\newtheorem{cor}[thm]{Corollary}
\newcommand{\thmref}[1]{Theorem~\ref{#1}}
\newcommand{\lemref}[1]{Lemma~\ref{#1}}

\newcommand{\propref}[1]{Proposition~\ref{#1}}

\theoremstyle{remark}
\newtheorem{rmk}{Remark}[section]

\begin{document}

\title[First ever negative Fourier coefficients]
{Oscillations and first-ever negative Fourier coefficients of symmetric square $L$-functions over sparse set}

\author{Srinivas  Kotyada and Lalit Vaishya}
\address{Srinivas (Kotyada)  Department of Mathematics, IISER Tirupati,  India}
\email{srinivas.kotyada@gmail.com}
\address{(Lalit Vaishya) Stat-Math Unit, Indian Statistical Institute, 7, S. J. S. Sansanwal Marg, New Delhi 110016, India}
\email{lalitvaishya@gmail.com; lalitvaishya$\_$24v@isid.ac.in}


\date{\today}
 
\maketitle

\begin{abstract}

 Let $sym^{2} f$ denote the symmetric square lift of  a Hecke eigenform $f \in S_{k}(\Gamma_{0}(N))$ with the $n^{\rm th}$-Fourier coefficients $ \lambda_{sym^{2}f}(n)$.  In this article, we prove an estimate for the first moment of the sequence $\{ \lambda_{sym^{2}f}(\mathcal{Q}(\underline{x}))\}_{\mathcal{Q} \in \mathcal{S}_{D}, \underline{x} \in \mathbb{Z}^{2}}$ where $\mathcal{S}_{D}$ denotes the set of in-equivalent reduced forms of the discriminant $D$. More precisely,  we establish an estimate for the following sum:
 \begin{equation*}
\begin{split}
S(sym^{2}f, D; X ) &=  \sideset{}{^{\flat }}\sum_{\substack{\mathcal{Q}(\underline{x}) \leq X \\ \underline{x} \in \mathbb{Z}^{2} ,~  \mathcal{Q} \in \mathcal{S}_{D} \\ \gcd(\mathcal{Q}(\underline{x}),N) =1  }} \lambda_{sym^{2}f}(\mathcal{Q}(\underline{x})),
\end{split}
\end{equation*}
 Moreover, we consider a question concerning the behavior of signs of the Fourier coefficients $\lambda_{sym^{2}f}(n),$ supported on the set of integers represented by  reduced forms of the discriminant $D$.   We determine the size of $n_{sym^{2}f, D}$  (see definition before \thmref{ExtMatKLSW}), in terms of the  conductor of the associated $L$-functions.
\end{abstract}

\section{Introduction and statements of the results}

The arithmetic behavior and distribution of values of arithmetical functions is the center of attraction in analytic number theory. The randomness in the behavior of these arithmetical functions leads to many interesting results. In this context, the arithmetical function  $n \mapsto \lambda_{f}(n)$,  where $\lambda_{f}(n)$ is the $n^{\rm th}$ Fourier coefficients of a Hecke eigenform $f$, is one of the interesting hosts which has attracted many mathematicians. There are several results such as non-vanishing, asymptotic of summatory function, and behavior of sign changes of the Hecke eigenvalues in the short interval.   The problems become more fascinating  if one considers the arithmetic behavior and distribution over the sparse set of positive integers. In this regard,  in a joint work with  Pandey \cite{LVMKP2021}, we established the estimates of the power moment and behavior of sign changes (in short intervals) of  $\lambda_{f}(n)$,  supported on the set of integers represented by a  reduced form of discriminant $D$ with certain assumption. In another work (see \cite{LalitManish}), we established an estimate for the first moment associated to a sequence of Hecke eigenvalues of a Hecke eigenform $f \in S_{k}(\Gamma_{0}(N))$, supported on the set of integers represented by a reduced form of discriminant $D$ with certain assumption. For more details, we refer to the introduction in  \cite{LalitManish}. In a similar vein, we prove an estimate for the first moment associated to $\lambda_{sym^{2}f}(n)$, the Fourier coefficients of symmetric square lift of a Hecke eigenform, supported on the set of integers represented by reduced forms of discriminant $D$.

Before stating our results, we fix some notations. 

\smallskip
Let $D$ denote negative discriminant and   $N$ be a square-free positive integer, and \linebreak $M_{k}(\Gamma_{0}(N), \chi)$ and $S_{k}(\Gamma_{0}(N), \chi)$ denote the ${\mathbb C}$- vector space of modular forms and cusp forms of integral weight $k$ for the congruence subgroup $\Gamma_{0}(N)$ with nebentypus $\chi$, respectively. Moreover, if the nebentypus character is trivial, we denote the space of modular forms and cusp forms by $M_{k}(\Gamma_{0}(N))$ and $S_{k}(\Gamma_{0}(N))$, respectively. A cusp form $f \in S_{k}(\Gamma_{0}(N))$ is said to be a Hecke eigenform (newform) if $f$ is a common eigenform for all the Hecke operators and the Atkin-Lehner $W$-operators. The space of   Hecke eigenforms (newforms) is denoted by $S^{new}_{k}(\Gamma_{0}(N))$. The Fourier coefficients of a Hecke eigenform satisfy the famous Ramanujan-Petersson bound which follows from the work of Deligne \cite{Deligne}. With these notations, let $f \in S_{k}(\Gamma_{0}(N))$ be a Hecke eigenform with the normalised $n^{\rm th}$-Hecke eigenvalue $\lambda_f(n)$. A cusp form $f$ is said to be a normalized form if $\lambda_f(1)=1$. The normalized Hecke eigenvalues $\lambda_{f}(n)'s$ are real, and the function $n \mapsto \lambda_{f}(n)$  is a multiplicative and satisfies the Hecke relation \cite[Eq. (6.83)]{Iwaniec}:
\begin{equation}\label{HeckeR}
\lambda_{f}(m)\lambda_{f}(n) = \sum_{d \vert \gcd(m,n)}  \lambda_{f}\left(\frac{mn}{d^2}\right)
\end{equation}
 for all positive integers $m$ and $n$ with $\gcd(mn, N)=1.$ From Deligne's estimate of Fourier coefficients of cusp form, we have 
\begin{equation} \label{lambda-coefficient-bound}
|\lambda_{f}(n)| \le d(n) \ll_{\epsilon} n^{\epsilon} 
\end{equation}
for  any arbitrary small  $\epsilon > 0$ where $d(n)$ denotes the number of positive divisors of $n$.  We define the Fourier coefficients $\lambda_{ sym^{2}f}(n)$ of a symmetric square lift $sym^{2}f$ of a normalised Hecke eigenform $f \in S_{k}(\Gamma_{0}(N))$ given by;
\begin{equation} \label{sym2-C}
\lambda_{ sym^{2}f}(n) = \displaystyle{\sum_{d^{2}m = n} \lambda_{f}(m^{2})}~ \quad \text { for all $n$ with }~   \gcd(n, N) =1. 
\end{equation}
From the construction, it is clear that the arithmetical function  $\lambda_{ sym^{2}f}(n)$ is real and multiplicative, and satisfies Deligne's estimate, and at prime $p$, we have $\lambda_{ sym^{2}f}(p) = \lambda_{ symf}(p^{2}) $.

\smallskip
Throughout the paper, $D$ denotes a negative discriminant and $\mathcal{Q}(\underline{x})$ denotes a primitive integral positive-definite binary quadratic form (reduced form) of fixed discriminant $D$  given by; $\mathcal{Q}(\underline{x}) = ax_1^2 + b x_1x_2 + c x_2^2,$ where $\underline{x} = (x_1,x_2)\in \mathbb{Z}^2$, $a, b, c \in \mathbb{Z}$ with $\gcd(a, b, c) = 1$ and discriminant $D = b^2-4ac$. For a given discriminant $D,$ let $\mathcal{S}_{D}$ denote the set of all in-equivalent reduced forms of discriminant $D$, i.e., 
$$
\mathcal{S}_{D} = \{ \mathcal{Q}: ~   \mathcal{Q} ~\text{is a reduced form of discriminant D}   \}/_\sim.
$$
 The set $\mathcal{S}_{D}$ forms a finite abelian  group under certain group operation. 
 Let $h(D)$ denote the class number, i.e., $ h(D)= \#\mathcal{S}_{D}.$ For details, we refer to \cite[chapter 2]{Cox}.

Our first result concerns about proving an upper bound estimate for the following sum defined as follows:
\begin{equation} \label{Upper-sym}
\begin{split}
S(sym^{2}f, D; X )  &: =   \sideset{}{^{\flat }}\sum_{\substack{\mathcal{Q}(\underline{x}) \leq X \\ \underline{x} \in \mathbb{Z}^{2} ,~  \mathcal{Q} \in \mathcal{S}_{D} \\ \gcd(\mathcal{Q}(\underline{x}),N) =1  }} \lambda_{sym^{2}f}(\mathcal{Q}(\underline{x})) \\
\end{split}
\end{equation}
where $\lambda_{ sym^{2}f}(n)$ denotes the Fourier coefficients of symmetric square lift of a Hecke eigenform $f$, and  $\mathcal{S}_{D}$ denotes the set of in-equivalent reduced forms of discriminant $D$. More precisely, we establish the following estimate.
 \begin{thm}\label{PropUpper}
Let $f \in S_{k}(\Gamma_{0}(N))$ be a normalised Hecke eigenform and $D$ be a discriminant. Then, 
for sufficiently large $X>0$ and any arbitrary small $\epsilon>0,$ we have  
\begin{equation*}
\begin{split}
S(sym^{2}f, D; X )  &
\ll_{\epsilon} (N^{2}k^{2}|D|^{\frac{3}{2}})^{\frac{1}{3}+\epsilon} X^{\frac{2}{3}+\epsilon}\\
\end{split}
\end{equation*}
where the implied constant depends only on $\epsilon$.
\end{thm}

A consequence of the above theorem, there is a result on the first negative Fourier coefficients of symmetric square lift of a Hecke eigenform if one can also have a lower bound estimate for the sum $S(sym^{2}f, D; X )$ under the certain assumption. Below, we give a brief description and state our second result.

\smallskip 
Let  $f \in S_{k}(\Gamma_{0}(N))$ be a normalized Hecke eigenform. The problem of finding the first negative Hecke eigenvalues of a Hecke eigenform $f$  is to establish an upper bound of $n_{f}$  (in terms of the analytic conductor) where  $n_{f}$ is the least positive integer   among all the positive integer $n$ satisfying 
$$
\lambda_{f}(n)<0 \quad {\rm and} \quad \gcd(n, N)=1. 
$$ 
This is a $GL_{2}$-analog of finding the least quadratic non-residue. In this direction, the bound for $n_{f}$ has been obtained subsequently in the following works; Kohnen and Sengupta \cite{Koh-Sen} ($n_{f} \ll (k^{2}N)^{26/60}$), Kowalski et al \cite{KLSW2010} ($n_{f} \ll (k^{2}N)^{9/20}$) and Matom$\ddot{\rm a}$ki \cite{KMatomaki} ($n_{f} \ll (k^{2}N)^{3/8}$).  The exponent in a result of  Matom$\ddot{\rm a}$ki (for the bound of $n_{f}$) is best known and merits her method. This bound of Matom$\ddot{\rm a}$ki is still far from the expected estimate as it is known that $n_{f} \ll (\log kN)^{2}$ under generalized Riemann hypothesis. 

To establish an upper bound for $n_{f}$, supported on the sparse set of integers, is a challenging task. In this direction,  adopting the method of Kowalski et al \cite{KLSW2010} and Matom$\ddot{\rm a}$ki \cite{KMatomaki},  in a joint work with Pandey  \cite{LalitManish}, we established an estimate for $n_{f, D}$, the first-ever sign change of  $\lambda_{f}(n)$ of a normalized Hecke eigenform $f$, supported on the set of integers represented by reduced form of discriminant $D$ with the class number $h(D)=1$ \cite[Theorem 1.2]{LalitManish} . An estimate for $n_{f,D}$ is given by
$$
n_{f, D} \ll_{\epsilon} (Nk^{2} |D|^{2})^{\frac{3}{4} + \epsilon}.
$$  
Here, we would like to point out that the above result is also valid if $h(D)>1$. Let $n_{sym^{2}f}$ denote the least positive integer among all the positive integers $n$ such that  $\lambda_{sym^{2}f}(n)<0$ and $\gcd(n, N)=1$, where $\lambda_{sym^{2}f}(n)$ denotes the $n^{\rm th}$ Fourier coefficients of the symmetric square lift $sym^{2}f$ of a normalised Hecke eigenform $f \in S_{k}(\Gamma_{0}(N))$. In the case of the symmetric square lift (a  $GL_3$-form) of a normalised Hecke eigenform $f$,  Lau et al \cite{YKLLW} established a subconvexity bound for $n_{sym^{2}f}$, i.e.,  $n_{sym^{2}f}\ll (k^{2}N^{2})^{40/113}$. As mentioned by Lau et al \cite{YKLLW}, the exponent $40/113~ (=0.35398)$ is smaller than $GL_{2}$-exponent $3/8$ of Matom$\ddot{\rm a}$ki.   They also provided a sufficient argument for it. 

Let $n_{sym^{2}f, D}$ denote the least integer among all $n$  such that $\lambda_{sym^{2}f}(n)<0$, $\gcd(n, N)=1$ and $n= \mathcal{Q}(\underline{x})$  for some $ \mathcal{Q} \in \mathcal{S}_{D}$,  $ \underline{x} \in \mathbb{Z}^{2}$.
 
 In our next result, we establish an upper bound estimate for $n_{sym^{2}f, D}$ following the idea of  Matom$\ddot{\rm a}$ki \cite{KMatomaki}.  We observe that it is not necessary that $n_{sym^{2}f} = n_{sym^{2}f, D}$ always,  as  $n_{sym^{2}f}$ can not always be represented by a reduced form of fixed discriminant $D$. Moreover, it is easy to see that  $n_{sym^{2}f} \le  n_{sym^{2}f, D}.$ Thus,  any bound for $ n_{sym^{2}f, D}$ would be sufficient to provide a bound for  $n_{sym^{2}f}$.  An estimate for $n_{sym^{2}f, D}$ is obtained in the following result.

\begin{thm}\label{ExtMatKLSW}
 Let  $f \in S_{k}(\Gamma_{0}(N))$ be a normalised Hecke eigenform and $D$ be a negative discriminant with class number $h(D)=1$ .  Then, for any arbitrarily small $\epsilon>0$, 
$$
n_{sym^{2}f, D} \ll_{\epsilon} (k^{2}N^{2} |D|^{\frac{3}{2}})^{\frac{61222}{100000} + \epsilon} {\sqrt{|D|}}^{5.51} 
$$
where the implied constant is absolute. Moreover, if $h(D)>1$, we have
$$
n_{sym^{2}f, D} \ll_{\epsilon} (k^{2}N^{2} |D|^{\frac{3}{2}})^{\frac{61222}{100000} + \epsilon} \left(\frac{\sqrt{|D|}}{h(D)} \right)^{\frac{9}{1.6334}}  
$$
 
\end{thm}

\begin{rmk} 

 It is important to note that the class number $h(D) =1 $ for the following fundamental discriminants $D = -3, -4, -8, -7, -11, -19, -43, -67, -163$ and non-fundamental discriminants $D =  -12, -16, -27, -28$ (see \cite[Remark 1.1]{Lalit})  
 \end{rmk}
 
We mentioned earlier that $n_{sym^{2}f} \le  n_{sym^{2}f, D}$. Hence, we have an estimate for $n_{sym^{2}f} $ in the following result. Taking $D=-4$, we have the following corollary. 
 \begin{cor}\label{Sym2sign}
 Let  $f \in S_{k}(\Gamma_{0}(N))$ be a normalised Hecke eigenform. Then, for any arbitrarily small $\epsilon>0$, 
$$
n_{sym^{2}f} \ll_{\epsilon} (k^{2}N^{2})^{0.61222 + \epsilon} 
$$
where the implied constant is absolute.
\end{cor}
\begin{rmk} 
The exponent $31/100$ for $n_{sym^{2}f}$ can be obtained following the argument of Matom$\ddot{\rm a}$ki \cite{KMatomaki}. It is better exponent in sub-convexity bound than the exponent ($40/113$) obtained by Lau et al  \cite{YKLLW}, and it  was mentioned in the work of  Matom$\ddot{\rm a}$ki \cite{KMatomaki} for $n_{sym^{2}f}$. Here in this work,  the associated L-function is a $GL_{3} \times GL_{2}$, so it is related to $GL_{3} \times GL_{2}$ exponent. This is the reason we are getting almost 2 times (in  \thmref{ExtMatKLSW}) the exponent  mentioned by Matom$\ddot{\rm a}$ki.
 \end{rmk}

\section{Key Ingredients}
In order to obtain our results, we need to study the sum $S(sym^{2}f, D; X )$ defined in \eqref{Upper-sym}. The sum $S(sym^{2}f, D; X )$ can be expressed in terms of known arithmetical functions given by 
\begin{equation}\label{PartialSum}
\begin{split}
S(sym^{2}f, D; X ) & 
= \sideset{}{^{\flat }}\sum_{n \le X \atop \gcd(n, N)=1}  (\lambda_{sym^{2}f}(n)) \sum_{ n =  \mathcal{Q}(\underline{x})  \atop \mathcal{Q} \in  \mathcal{S}_{D}~ \underline{x} \in \mathbb{Z}^{2}} 1 =  \sum_{n \le x } \lambda_{sym^{2}f}(n) r_{D}(n)\\
\end{split}
\end{equation}
where symbol $\flat $ means that the sum runs over all square-free positive integers and $r_{D}(n)$ denotes the number of representations of $n$ by reduced forms of discriminant $D$. 
We define the generating function  $\theta_{D}(\tau)$ associated to  $r_{D}(n)$ as follows:  
\begin{equation*}
\begin{split}
\theta_{D}(\tau) & := \quad  \displaystyle{\sum_{n = 0}^{\infty} r_{D}(n) q^{n}}, \qquad   \quad ~~~ q = e^{2 \pi i \tau}. \\
\end{split}
\end{equation*} 
In the above expression, the constant term $r_{D}(0)$ is precisely the class number $h(D)$. The generating function $\theta_{D}(\tau) \in M_{1}(\Gamma_{0}(|D|), \chi_{D}) $ (see \cite[Theorem 10.9]{Iwaniec}), where $\chi_{D}$ is the Dirichlet character of  modulus $|D|$ which  is given by Jacobi symbol $\chi_{D}(d) := \left(\frac{D}{d} \right)$.  From  Weil's bound,  it is well-known that  $r_{D}(n) \ll n^{\epsilon}$ for any arbitrarily small  $\epsilon > 0.$  In the literature, the quantity $r_{D}(n)$ is also known as character sum.
A  formula for $r_{D}(n),$ the number of representations of $n$ by reduced forms of discriminant $D,$ is given by \cite[section 11.2]{Iwaniec} 
\begin{equation} \label{r-quad}
r_{D}(n) = w_{D} \sum_{d \vert n} \chi_{D}(d), \   \
{\rm ~where~} \ \ \ \ \ 
w_{D} =
\begin{cases} 
6 {\rm ~~~~~ if~~~~~~} D = -3,\\
4 {\rm ~~~~~ if~~~~~~} D = -4,\\
2 {\rm ~~~~~ if~~~~~~} D < -4. \end{cases} 
\end{equation}
Here,  $w_{D}$ denotes the size of the group of units of the imaginary quadratic field $\mathbb{Q}(\sqrt{D})$.  The formula for $r_{D}(n)$ depends only on the  discriminant $D$ but does not depend on the choice of  reduced forms appearing in the class group corresponding to discriminant $D.$  Let us write \linebreak  $r_{D}(n) = w_{D} r^{*}_{D}(n)$ where $r^{*}_{D}(n) = \displaystyle{\sum_{d \vert n}}~ \chi_{D}(d)$. Then, the sum in \eqref{PartialSum} can be written as follows:
\begin{equation}\label{Sum-PIBQF}
\begin{split}
S(sym^{2}f, D; X ) 
& =  \sideset{}{^{\flat }} \sum_{n \le X \atop \gcd(n,N) =1 } \lambda_{sym^{2}f}(n) r_{D}(n) = w_{D} \sideset{}{^{\flat }} \sum_{n \le X \atop \gcd(n,N) =1 } \lambda_{sym^{2}f}(n) r^{*}_{D}(n).
\end{split}
\end{equation}
We define the following Dirichlet series associated to $\lambda_{sym^{2}f}(n)$ and  $r^{*}_{D}(n)$ given by  
\begin{equation}\label{R-L-function}
\begin{split}
L(sym^{2}f, D; s) &=  \sideset{}{^{\flat }} \sum_{n \ge 1 \atop \gcd(n,N) =1 } \frac{\lambda_{sym^{2}f}(n) r^{*}_{D}(n)}{n^{s}} 
\end{split}
\end{equation}
The Dirichlet series  $L(sym^{2}f, D; s)$ converges absolutely and uniformly for $\Re(s) >1.$   To obtain an upper bound for the sum $S(sym^{2}f, D; X ),$ we first decompose $L(sym^{2}f, D; s )$ in terms of known $L$-functions. Using the analytic properties of the associated $L$-functions, we establish our estimate. 
Let $f(\tau) = \displaystyle{\sum_{n=1}^\infty \lambda_f(n)n^{\frac{k-1}{2}}q^{n}} \in S_{k}(\Gamma_{0}(N))$ be a normalised Hecke eigenform and $\psi_{N}$ denote the principal character of  modulus $N$.
We define the symmetric square $L$-function associated to $f$ as follows:
\begin{equation*}
\begin{split}
L(sym^{2}f, s) &= \sum_{n \ge 1} \frac{\lambda_{sym^{2}f}(n)}{n^{s}} =  \prod_{p}\left(1-\frac{\lambda_{f}(p^{2})}{p^{s} } +\frac{\psi_{N}(p)\lambda_{f}(p^{2})}{p^{2s} }+ \frac{\psi_{N}(p)}{p^{3s}}\right)^{-1} 
\end{split}
\end{equation*}
Similarly, we define $L(sym^{2}f \otimes \chi_{D}, s)$,  the twist of  
symmetric square $L$- function $L(sym^{2}f, s) $ by a Dirichlet character $\chi_{D}$ of modulus $D$, given by
\begin{equation*}
\begin{split}
L(sym^{2}f \otimes \chi_{D}, s) &= \sum_{n \ge 1} \frac{\lambda_{sym^{2}f}(n) \chi_{D}(n)}{n^{s}}. 
\end{split}
\end{equation*}
These Dirichlet series converge absolutely for $\Re(s) > 1$ and do not vanish in the region $\Re(s) \ge 1$. Moreover, we define the associated completed $L$-functions given by
\begin{equation*}
\begin{split}
\Lambda( sym^{2}f,s) &:= N^{s}\Gamma \left( \frac{s+1}{2}\right)  \Gamma \left( \frac{s+k-1}{2}\right) \Gamma \left( \frac{s+k}{2}\right)  L(sym^{2}f,s)  \\
\end{split}
\end{equation*}
and
\begin{equation*}
\begin{split}
\Lambda( sym^{2}f\otimes \chi_{D},s) &:= (N^{2}|D|^{3})^{\frac{s}{2}} \Gamma  \! \left( \frac{s+1}{2}\right)  \Gamma \! \left( \frac{s+k-1}{2}\right) \Gamma \! \left( \frac{s+k}{2}\right) L(sym^{2}f \otimes \chi_{D},s).
\end{split}
\end{equation*}
These completed $L$-functions are analytically continued to the whole $\mathbb{C}$-plane and satisfy a nice functional equation $(s \rightarrow 1-s)$. For details, we refer to \cite[Theorem 1]{GS1975}.  The decomposition of  $L(sym^{2}f, D; s)$ is given in the following Lemma.
\begin{lem}\label{Decomp}
For $\Re(s) >1 $, we have 
\begin{equation} \label{U(s)}
\begin{split}
L(sym^{2}f, D; s) & = L( sym^{2}f,s) L( sym^{2}f \otimes {\chi_{D}},s) G(s), \\ 
\end{split}
\end{equation}
where $L( sym^{2}f,s)$ and  $L( sym^{2}f \otimes {\chi_{D}},s)$ are the symmetric square $L$-function and the  twisted symmetric square  $L$-function, respectively, and the Euler product for $G(s)$ is given by
\begin{equation*}
\begin{split}
\\
G(s) &=  \prod_{p \mid N|D|^{2}}\left( 1- \frac{ \lambda_{f}(p)}{p^{s}}\right)    \prod_{p \nmid N|D|^{2}}  \bigg( 
  1+ \frac{2 - 2 \lambda^{2}_{f}(p)- \lambda^{2}_{f}(p) \chi_{D}(P) }{p^{2s} }  \\
& +  \frac{ \lambda_{f}(p)(1+\chi_{D}(p)) (1+ \lambda^{2}_{f}(p) \chi_{D}(P)) }{p^{3s}} +  \frac{ 1- 2\lambda^{2}_{f}(p)(1+ \chi_{D}(P)) }{p^{4s} } + \frac{\lambda_{f}(p)(1+ \chi_{D}(P))}{p^{5s}}
  \bigg).
\end{split}
\end{equation*}
It converges absolutely and uniformly in the half plane $\Re(s) > \frac{1}{2}$ and non-zero for $\Re(s) = 1$.
\end{lem}

\noindent
\textbf{Proof:} ~
We prove the above lemma following similar arguments as in the proof of \cite[Lemma 2.1]{Lalit}.
We know that $\lambda_{sym^{2}f}(n)$ and  $r^{*}_{D}(n)$ are multiplicative functions. So, for $\Re(s) >1,$ we have 
\begin{equation*}
\begin{split}
L(sym^{2}f, D; s ) &=  \!\!\!\! \sideset{}{^{\flat }} \sum_{n \ge 1 \atop \gcd(n,N) =1 } \frac{\lambda_{sym^{2}f}(n) r^{*}_{D}(n)}{n^{s}} = \prod_{p \nmid N}\left( 1+ \frac{\mu^{2}(p) \lambda_{sym^{2}f}(p) r^{*}_{D}(p) }{p^{s}}\right) \\
& =  \prod_{p \nmid N}\left( 1+ \frac{ \lambda_{sym^{2}f}(p) (1+ \chi_{D}(p)) }{p^{s}}\right)= L( sym^{2}f,s) L( sym^{2}f \otimes \chi_{D},s )G(s), \\
\end{split}
\end{equation*}
where the Euler product for $G(s)$ is given in \lemref{Decomp}.


\begin{lem} \cite[Chapter 5]{HIwaniec}
Let  $f \in S_{k}(\Gamma_{0}(N))$ be a normalised Hecke eigenform. The hybrid convexity bound of the symmetric square  $L$-function and the  twisted symmetric square  $L$-function is given by:
\begin{equation*}
\begin{split}
L( sym^{2}f,s) &\ll  \left( N^{2} k^{2} (|s|+3)^{3}\right)^{\frac{1}{2}(1-\sigma)+\epsilon} ~ \\ 
\end{split}
\end{equation*}
and
\begin{equation*}
\begin{split}
L( sym^{2}f \otimes {\chi_{D}}, s) &\ll \left( N^{2}  k^{2} |D|^{3} (|s|+3)^{3}\right)^{\frac{1}{2}(1-\sigma)+\epsilon} ~ \\ 
\end{split}
\end{equation*}
where $\sigma = \Re(s)$ and $\frac{1}{2} \le \sigma \le 1$.
\end{lem}

\begin{lem}\label{Mean-value} \cite[Theorem 5]{APerelli}
Let $L(F,s)$ be an $L$-function of degree $m \ge 2$. Then we have 
\begin{equation}
\begin{split}
\int_{T}^{2T} \left|L \left(F, \sigma+it \right)\right|^{2} dt &\ll_{ \epsilon} (Q_{F}(1+|t|)^{m})^{1-\sigma +\epsilon}
\end{split}
\end{equation}
uniformly for $\frac{1}{2} \le \sigma \le 1$ and $|t| \ge 1$, and $Q_{F}$ denotes the analytic conductor of $F$. 
\end{lem}

\subsection{ \textbf{Mean value of a multiplicative function over the integers represented by a reduced form of discriminant $D$}} 
 \begin{lem}\cite[Propostion 2.3]{LalitManish}\label{MVEBQF}
Let $\eta$ be a fixed positive integer and $r^{*}_{D}(n)$ is given by; $r^{*}_{D}(n) = \displaystyle{\sum_{d \vert n}}~ \chi_{D}(d).$
Then there exist an absolute constant $C = C(\eta)$ such that 
\begin{equation}\label{MFEST}
\begin{split}
\mathcal{E}_{\eta}( X )  = \!\!\!\! \sum_{n \le X \atop \gcd(n,N) =1} \!\!\!\! \mu^{2}(n) \eta^{\omega(n)} r^{*}_{D}(n)
= \frac{P(1) L(1, \chi_{D})^{\eta}}{\Gamma(\eta)} X (\log X)^{\eta-1} \left(  1+ O_{\eta} \left( \frac{L_{N}^{2 e \eta +2}}{\sqrt{\log X}} \right)\right)
\end{split}
\end{equation}  
 uniformly for $N \ge 1$ and $X \ge \exp (C L_{N}^{2 e \eta +2} )$ where $ L_{N} = \log (\omega(N)+3)$ and
 \begin{equation*}
\begin{split}
P(s)  & =  \prod_{p \mid N} \left[\left( 1- \frac{1}{p^{s}}\right) \left( 1- \frac{\chi_{D}(p)}{p^{s}}\right)\right]^{\eta}  \prod_{p \nmid N} \left [ \left( 1- \frac{1}{p^{s}}\right)^{\eta} \!\! \left( 1- \frac{\chi_{D}(p)}{p^{s}}\right)^{\eta}  \!\!
\left(1+ \frac{\eta (1+\chi_{D}(p))}{p^{s}}\right) \right]. 
\end{split}
\end{equation*}
The Euler product for $P(s)$ converges absolutely and uniformly for $\Re(s) > 1/2$ and non-zero at $s=1.$ This result has been proved in \cite[Propostion 2.3]{LalitManish} whenever $h(D)=1$ but the result is also valid when $h(D)>1$.  
\end{lem}  

\section{Proof of Results}
\noindent
\subsection{General Philosophy to obtain an upper bound estimate for the first moment:} \quad  Let $1 \le P \le \frac{X}{2}.$
In order to obtain an upper bound for the sum $S(sym^{2}f, D; X )$ given in \eqref{Sum-PIBQF}, we introduce a smooth compactly supported function $w(x)$ satisfying: $w(x) =1$ for $x \in [2P, X],$ $w(x) = 0$ for $x<P$ and $x> X+P,$ and $w^{(r)}(x) \ll_{r} P^{-r}$ for all $r\ge 0.$ In general, for any arithmetical function $n \mapsto f(n), $ we have 
\begin{equation}\label{UpperEstimate}
\sum_{n \le X} f(n) =  \frac{1}{2 \pi i} \int_{(b)} \tilde w(s) \left(\sum_{n\ge 1} \frac{f(n)}{n^{s}}\right)  ds +  O\left(\sum_{n< 2P} |f(n)| \right) + O\left(\sum_{X< n< X+ P} |f(n)| \right). 
\end{equation}
where $b$ is a positive real number larger than the  abscissa of absolute convergence of \linebreak $\displaystyle{\sum_{n\ge 1}{f(n)}{n^{-s}}}$ and $\tilde w(s)$ denotes the Mellin's transform of  $ w(t)$ is given by following integral; 
$$\tilde w(s) = \int_{0}^{\infty} w(x) x^{s} \frac{dx}{x}.$$
We observe that due to integration by parts, 
\begin{equation}\label{FourierW}
\tilde w(s) =  \frac{1}{s(s+1)\cdots(s+t-1)}\int_{0}^{\infty} w^{(m)}(x) x^{s+m-1} dx \ll \frac{P}{X^{1-\sigma}}  \left(\frac{X}{|s|P}\right)^{m}
\end{equation}
for any $m\ge 0,$ where $\sigma = \Re(s).$ For details, we refer to \cite[Section 3]{LalitManish}.

\smallskip
\noindent
\textbf{Proof of \thmref{PropUpper}}
 From Equation \eqref{UpperEstimate} with $f(n) = \lambda_{sym^{2}f}(n)  r^{*}_{D}(n),$ we have 
\begin{equation}\label{GenSum}
\begin{split}
\sideset{}{^{\flat }} \sum_{n \le X \atop \gcd(n,N)=1} \!\!\!\! \lambda_{sym^{2}f}(n)  r^{*}_{D}(n) & = 
\frac{1}{2 \pi i} \int_{(1+\epsilon)} \!\! \tilde w(s) L(sym^{2}f, D; s )  ds \\
& ~ +  O \!\!\left(\sideset{}{^{\flat }} \sum_{n < 2P \atop \gcd(n,N) =1 } \!\!\!\!\! |\lambda_{sym^{2}f}(n)  r^{*}_{D}(n)|  \!\! \right) + O\!\! \left(\sideset{}{^{\flat }} \sum_{X< n< X+ P \atop \gcd(n,N) =1 }  \!\!\!\!\! |\lambda_{sym^{2}f}(n)  r^{*}_{D}(n)|  \!\! \right). 
\end{split}
\end{equation}
Since   $ \lambda_{sym^{2}f}(n) \ll n^{\epsilon}$   and  $r^{*}_{D}(n)  \ll n^{\epsilon}$ for any $\epsilon>0$. So, we have  
\begin{equation}\label{Errorbound}
\begin{split}
  O\left(\sideset{}{^{\flat }} \sum_{n < 2P \atop \gcd(n,N) =1 } |\lambda_{sym^{2}f}(n)  r^{*}_{D}(n)| \right) + O\left(\sideset{}{^{\flat }} \sum_{X< n< X+ P \atop \gcd(n,N) =1 } |\lambda_{sym^{2}f}(n)  r^{*}_{D}(n)| \right)  \ll P^{1+\epsilon}.
\end{split}
\end{equation}
Moreover, $\tilde w(s)  \ll \frac{P}{X^{1-\sigma}}  \left(\frac{X}{|s|P}\right)^{m}$
for any $m \ge 0,$  the contribution for the integral (appearing in \eqref{GenSum}) over $|s| \ge T = \frac{X^{1+\epsilon}}{P}$, is negligibly small, i.e., $O(X^{-A})$ for any large $A>0$ if one chooses sufficiently large $m>0.$ More precisely,
\begin{equation*}
\frac{1}{2 \pi i} \int_{(1+\epsilon)} \!\! \tilde w(s) L(sym^{2}f, D; s )  ds = \frac{1}{2 \pi i} \int_{1+\epsilon-iT}^{1+\epsilon+iT} \tilde w(s) L(sym^{2}f, D; s )  ds  + O(X^{-A}). 
 \end{equation*}
Now, we shift the line of integration from $\Re(s) = 1+\epsilon$ to $\Re(s) =\sigma_{0}:= \frac{1}{2}+\epsilon$ and apply Cauchy's residue theorem and standard argument to get 
\begin{equation}\label{RHSInt}
\frac{1}{2 \pi i} \int_{(1+\epsilon)} \!\! \tilde w(s) L(sym^{2}f, D; s )  ds =  \frac{1}{2 \pi i} \int_{\sigma_{0}-iT}^{\sigma_{0}+iT} \tilde w(s) L(sym^{2}f, D; s )  ds + O(X^{-A})  . 
 \end{equation}
We substitute the decomposition $L(sym^{2}f, D; s)  = L( sym^{2}f,s) L( sym^{2}f \otimes {\chi_{D}},s) G(s)$  from \lemref{Decomp},  bound for $\tilde w(s)$ when $m=1$ (given in \eqref{FourierW}),  and use absolute convergence of $G(s)$ in the region $\Re(s)> 1/2$ to get
\begin{equation*}
\begin{split}
 & \frac{1}{2 \pi i} \int_{\sigma_{0}-iT}^{\sigma_{0}+iT} \tilde w(s) L(sym^{2}f, D; s )  ds  
   \ll \int_{-T}^{T} |\tilde w(\sigma_{0}+it)| |L(sym^{2}f, D; \sigma_{0}+it )| dt  \\
  & \qquad  \ll  \int_{0}^{T} \frac{X^{\sigma_{0}}}{|\sigma_{0}+it|} |L(sym^{2}f, D; \sigma_{0}+it )| dt  \\
  & \qquad \ll \left(\int_{0}^{1} + \int_{1}^{T}\right) \frac{X^{\sigma_{0}}}{|\sigma_{0}+it|} |L(sym^{2}f, \sigma_{0}+it ) L( sym^{2}f \otimes \chi_{D}, \sigma_{0}+it  )| dt  \\
  & \qquad \ll (N^{2}k^{2} |D|^{\frac{3}{2}} X)^{\sigma_{0}} + X^{\sigma_{0}} \!\!\! \int_{1}^{T}  \frac{|L( sym^{2}f, \sigma_{0}+it  ) L( sym^{2}f \otimes \chi_{D}, \sigma_{0}+it  )|}{t} dt.\\
  \end{split}
 \end{equation*}
We apply the dyadic division method and then Cauchy-Schwartz inequality to get   
\begin{equation*}
\begin{split}
 & \int_{1}^{T}  \frac{|L( sym^{2}f, \sigma_{0}+it  ) L( sym^{2}f \otimes \chi_{D}, \sigma_{0}+it  )|}{t} dt  \\
 &   \ll \log T  \underset{1 \le T_{1} \le T}{\rm max} \!\! \left\{ \frac{1}{T_{1}} \!\!  \left( \int_{\frac{T_{1}}{2}}^{T_{1}} \!\!\!\!  |L( sym^{2}f, \sigma_{0}+it)|^{2} dt \right)^{\frac{1}{2}} \!\! \left(\int_{\frac{T_{1}}{2}}^{T_{1}} \!\!\!\! |L( sym^{2}f \otimes \chi_{D}, \sigma_{0}+it)|^{2} dt \right)^{\frac{1}{2}}\right\}   \\
 & \ll \log T  \underset{1 \le T_{1} \le T}{\rm max} \left(\frac{1}{T_{1}} (N^{2}k^{2} T_{1}^{3})^{\frac{1}{4}+\epsilon} (N^{2}k^{2}  |D|^{3} T_{1}^{3})^{\frac{1}{4}+\epsilon}\right)
   \ll \log T (N^{2}k^{2}  |D|^{\frac{3}{2}})^{\sigma_{0}} T^{1/2}.
  \end{split}
 \end{equation*}
Now, we combine all the estimates in \eqref{RHSInt} to get 
\begin{equation}\label{halfInt}
\begin{split}
\frac{1}{2 \pi i} \int_{(1+\epsilon)} \!\!\!\! \tilde w(s) L(sym^{2}f, D; s )  ds 
& =  O\left((N^{2}k^{2}  |D|^{\frac{3}{2}} X)^{\sigma_{0}} \right)  + O\left( (N^{2}k^{2} |D|^{\frac{3}{2}} )^{\sigma_{0}} X^{1+\epsilon} T^{\frac{1}{2}} \right) \\ 
& =O\left( (N^{2}k^{2} |D|^{\frac{3}{2}} )^{ \frac{1}{2}+\epsilon} X^{ \frac{1}{2}+\epsilon} T^{\frac{1}{2}} \right).
\end{split} 
\end{equation}
Here, $\sigma_{0}= \frac{1}{2}+\epsilon$. From equations \eqref{GenSum}, \eqref{Errorbound}, \eqref{RHSInt} and \eqref{halfInt} with  $T = \frac{X^{1+\epsilon}}{P}$, we have 
 \begin{equation*}
\begin{split}
 \sideset{}{^{\flat }} \sum_{n \le X \atop \gcd(n,N)=1} \!\!\!\! \lambda_{sym^{2}f}(n)  r^{*}_{D}(n) 
 & =  O(P^{1+\epsilon})+ O\left( (N^{2}k^{2} |D|^{\frac{3}{2}} )^{\frac{1}{2}+\epsilon} X^{1+\epsilon} P^{-\frac{1}{2}} \right).
\end{split}
\end{equation*}
Now, we choose $P= (N^{2}k^{2}  |D|^{\frac{3}{2}})^{\frac{1}{3}+\epsilon} X^{\frac{2}{3}+\epsilon}$ to get 
 \begin{equation*}
\begin{split}
 \sideset{}{^{\flat }} \sum_{n \le X \atop \gcd(n,N)=1} \!\!\!\! \lambda_{sym^{2}f}(n)  r^{*}_{D}(n) 
 & =  O\left( (N^{2}k^{2} |D|^{\frac{3}{2}} )^{\frac{1}{3}+\epsilon} X^{\frac{2}{3}+\epsilon}  \right). 
\end{split}
\end{equation*}

\subsection{A positive lower bound for $S(sym^{2}f, D; Y^{u} )$:}
Let $Y$ be a real number such that $Y \le X/2$ and $\lambda_{sym^{2}f}(n) \ge 0$ for all $n < Y$ with $\gcd(n,  N)=1$ and $n= \mathcal{Q}(\underline{x})$  for some $ \mathcal{Q} \in \mathcal{S}_{D}$ and  $ \underline{x} \in \mathbb{Z}^{2}$.
We obtain an estimate for $n_{sym^{2}f, D}$ by comparing positive lower bound and upper bound estimates for the sum $S(sym^{2}f, D; Y^{u} )$ (defined in \eqref{Upper-sym}), for some $u> 1$. An upper bound is comparatively easy and it is obtained in \thmref{PropUpper}, using the $L$-function techniques. In order to obtain a lower bound for  $S(sym^{2}f, D; Y^{u} )$, for some $u >  1$, one needs to develop new tools and invoke the following identity (as mentioned in \cite{YKLLW}):
$$
\lambda_{sym^{2}f}(p^{\nu}) = \frac{\sin ((\nu+2)\theta_{p}) \sin ((\nu+1)\theta_{p})}{\sin (\theta_{p})\sin (2\theta_{p})} \quad {\rm for ~ all~ } p \nmid N ~{\rm and }~ \nu \ge 1.
$$
Moreover, one  needs to exclude those primes for which $\lambda_{sym^{2}f}(p^{\nu}) = 0$ for all primes $p$ and integers $ \nu \ge 1$ with $p^{\nu} \le Y$. For $\nu \le 4$,  Lau et al \cite{YKLLW} mentioned that there are only a few primes because it is equivalent to enumerate primes $p$ such that $\lambda_{f}(p) = \alpha$ for some non-zero real algebraic number $\alpha$. More precisely,  they proved that
$$\#\{p: \lambda_{f}(p) = \alpha\}  \ll \log(kN)$$
for all Hecke eigenform $f$ and  non-zero real algebraic number $\alpha$  \cite[Lemma 2.4]{YKLLW}. They utilized these arguments to establish an estimate for $n_{sym^{2}f}$.
 We define the following set
$$
\mathcal{P}^{*}_{f} :=\underset{1 \le \nu \le \log Y} \cup \{ p: \lambda_{f}(p) = 2 \cos(\pi/(\nu+2))\}
$$
Following the argument of  Lau et al \cite{YKLLW}, we have 
$$
\#\mathcal{P}^{*}_{f} \ll (\log (kN))^{2}
$$
for all Hecke eigenform $f \in S_{k}(\Gamma_{0}(N)) $ whenever $\log Y \ll \log (kN)$. Let us consider the set $N_{f} = \{p \mid N\} \cup \{p \in \mathcal{P}^{*}_{f} \}$. Then, there are few primes in $N_{f}$ as $\# \{ p \in N_{f}  \} \ll  (\log (kN))^{2}.$  Following the idea of   Kowalski et al \cite{KLSW2010} and  Lau et al \cite{YKLLW}, we make use of   an auxiliary multiplicative function $h_{Y}$, supported on the set of square-free integers, defined by
\begin{equation}\label{Auxfun}
\begin{split}
&h_{Y}(p) = 
\begin{cases}
\beta \left( \frac{\log p}{\log Y} \right),\qquad {\rm if} \quad  p< Y \quad {\rm and } \quad p\not\in N_{f}, \\
-1, \qquad \qquad \quad {\rm if} \quad  p \ge Y \quad {\rm or } \quad p\not\in N_{f}, \\
0, \qquad \qquad \qquad {\rm if} \quad   p\in N_{f} \quad {\rm or} \quad  p>Y , \\
\end{cases} \\
\end{split}
\end{equation}
where the function $\beta: [0, \infty) \rightarrow \mathbb{R}$ is given by
\begin{equation}\label{BETA}
\begin{split}
 \beta(t) =
\begin{cases}
3, \qquad  \qquad \qquad \qquad  {\rm if } \qquad    t=0,\\
3-4\sin^{2}\left(\frac{\pi}{m+2}\right) , \quad ~\quad {\rm if } \quad \frac{1}{m+1} \le  t < \frac{1}{m}, ~ ~m \in \mathbb{N},\\
 -1, \qquad \quad \quad \quad \quad \quad    {\rm if} \qquad t \ge 1.
\end{cases}
\end{split}
\end{equation}
 We make use of the auxiliary function $h_{Y}$  to obtain a positive lower bound for the sum \linebreak  $S(sym^{2}f, D; Y^{u} )$ for some $u > 1$. 
 
 By definition of $Y$,  $\lambda_{sym^{2}f}(n) \ge 0$ for every square-free positive integer $n < Y,$ with $\gcd(n,N) =1$ such that  $n= \mathcal{Q}(\underline{x})$  for some $ \mathcal{Q} \in \mathcal{S}_{D}$ and  $ \underline{x} \in \mathbb{Z}^{2}$. Then, the Hecke relation together with the assumption  $\lambda_{sym^{2}f}(p) \ge 0$ for every prime $p$ and $m \ge 1$ with $p \le Y^{1/m}$ implies that $ \lambda_{sym^{2}f}(p) \ge 3 -4 \sin^{2}\left(\frac{\pi}{m+2}\right).$ This shows that $\lambda_{sym^{2}f}(p) \ge h_{Y}(p),$ for all primes $p < Y.$ Multiplicativity of these arithmetical functions allows us to get $\lambda_{sym^{2}f}(n) \ge h_{Y}(n),$ for all square-free integers $n < Y$. Similar to \cite{KLSW2010}, we  observe that for each $u \ge u_{0}= 1.6334$ ($ u_{0}= 1/0.61222$), 
\begin{equation}\label{LBound}
\begin{split}
S(sym^{2}f, D; Y^{u}) & = \sideset{}{^{\flat }} \sum_{n \le Y^{u} \atop \gcd(n,N) =1 } \lambda_{sym^{2}f}(n) r^{*}_{D}(n) \ge \sideset{}{^{\flat }} \sum_{n \le Y^{u} \atop \gcd(n,N) =1 } h_{Y}(n) r^{*}_{D}(n)>0.
\end{split}
\end{equation}
Precisely, to see \eqref{LBound}, let $g_{Q}$ be a multiplicative function defined by convolution identity $\lambda^{*}_{sym^{2}f}(n)  =  (g_{Q} * h^{*}_{Y})(n)$,
where $\lambda^{*}_{sym^{2}f}(n) = \lambda_{sym^{2}f}(n)r^{*}_{D}(n)$ and  $h^{*}_{Y}(n) = h_{Y}(n)r^{*}_{D}(n). $ This gives that  $g_{Q}(p) =\lambda_{sym^{2}f}(p)r^{*}_{D}(p)- h_{Y}(p)r^{*}_{D}(p).$ If $p$ is not represented a reduced form of discriminant $D,$  then $g_{Q}(p) =0.$ If $p$ is represented by some reduced form of discriminant $D,$ then $g_{Q}(p) > 0$ for  such $p$ not dividing $N_{f}$ with $p<Y$, by definition of  $ h_{Y}(p).$ For $p \ge Y,$ with $p \nmid N,$ we have $h_{Y}(p)r^{*}_{D}(p) = -2$ and $-1\le \lambda_{sym^{2}f}(p) \le 3.$ Hence, $g_{Q}(p) = 2a+2 \ge 0$ as  $-1\le a =\lambda_{sym^{2}f}(p) \le 3.$ Thus, we have $g_{Q}(n) \ge 0$ for all square-free positive integers $n.$ Now, we see that for each $u \ge u_{0}$,
\begin{equation*}
\begin{split}
S(sym^{2}f, D; Y^{u}) & = \sideset{}{^{\flat }} \sum_{n \le Y^{u} \atop \gcd(n,N) =1 } \lambda_{sym^{2}f}(n) r^{*}_{D}(n) =  \sideset{}{^{\flat }} \sum_{n \le Y^{u} \atop \gcd(n,N) =1 } \sum_{d \mid n} g(d) h_{Y}(n/d) r^{*}_{D}(n/d)\\
& =  \sideset{}{^{\flat }} \sum_{d \le Y^{u} \atop \gcd(d,N) =1 }  g(d)  \sideset{}{^{\flat }} \sum_{m \le Y^{u}/d \atop \gcd(m,N) =1 }   h_{Y}(m) r^{*}_{D}(m) \ge \sideset{}{^{\flat }} \sum_{n \le Y^{u} \atop \gcd(n,N) =1 }   h_{Y}(n) r^{*}_{D}(n)>0
\end{split}
\end{equation*} 
as $g(1) =1 $ and $g(n )\ge 0$ for each square-free positive integer $n.$ Thus, a positive lower bound of $S(sym^{2}f, D; Y^{u_{0}})$ is given by  $\displaystyle{\sideset{}{^{\flat }} \sum_{n \le Y^{u_{0}} \atop \gcd(n, N) =1 }   h_{Y}(n) r^{*}_{D}(n)}$ for some $u_{0}$. For $u \ge u_{0}$, an estimate for $\displaystyle{\sideset{}{^{\flat }} \sum_{n \le Y^{u} \atop \gcd(n, N) =1 }   h_{Y}(n) r^{*}_{D}(n)}$ is given in the following proposition.

Let $\beta$ be a step function (defined in \eqref{BETA}) given by;  $\beta: [0, \infty) \rightarrow \mathbb{R}$ with $\beta(t) = \beta_{k}$ when $t \in [x_{k}, x_{k+1}]$ for each $k= 0,1,2, \cdots K$.  Here $K > 0$ and $[x_{k}, x_{k+1}]$ is partition of $[0, \infty)$. In our case $x_{0}=0$, $x_{k} = \frac{1}{k}$,  $\beta_{0} =2$ and $\beta_{k} = 3- 4 \sin^{2} \left(\frac{\pi}{k+2}\right)$.
  \begin{prop}\label{PropLower}
Let $U \ge 1$ be a real number and $D$ be a discriminant with $h(D)=1$, and  $h_{Y}(n)$ (defined in \eqref{Auxfun}) and $\beta(t)$ be as above with $\beta_{0}>0$. Let $N \le X^{U}$ be a positive integer. Then, we have
\begin{equation*}
\begin{split}
 \sum_{n \le Y^{u} \atop \gcd(n,N) =1 }   h_{Y}(n) r^{*}_{D}(n) &= (\sigma(u) + o_{\beta, U}(1)) \frac{P(1) L(1, \chi_{D})^{\beta_{0}}}{\Gamma(\beta_{0})} (\log Y^{u})^{\beta_{0}-1} Y^{u} 
\end{split}
\end{equation*}
uniformly for $u \in \left[\frac{1}{U}, U \right]$ where $P(s)$ is given in \lemref{MVEBQF} and 
\begin{equation*}
\begin{split}
\sigma(u)  &= u^{\beta_{0}-1} + \sum_{j=1}^{\infty} \frac{(-1)^{j}}{j}  I_{j}(u)  \\
\end{split}
\end{equation*}
with
\begin{equation*}
\begin{split}
 I_{j}(u) &= \int_{\Delta_{j}} (u- t_{1}- t_{2}- \cdots- t_{j})^{\beta_{0}-1} \prod_{i=1}^{j} (\beta_{0}- \beta(t))  \frac{dt_{1} dt_{2} \cdots dt_{j}}{ {t_{1}t_{2} \cdots t_{j}}} \\
\end{split}
\end{equation*}
and
\begin{equation*}
\begin{split}
  & \Delta_{j} = \{ ( t_{1}, t_{2}, \cdots, t_{j}) \in [0, \infty) \mid  |t_{1}+ t_{2}+ \cdots+t_{j}| \le u \}.
\end{split}
\end{equation*}
This result also hold when $h(D)>1$.
\end{prop}
The proof follows exactly the same argument as in the proof of \cite[Proposition 2.4]{LalitManish}. In \cite[Proposition 2.4]{LalitManish},  the above proposition (with a different function $\alpha(t)$) has been proved assuming that the class number $h(D)=1$ but the result is also valid when $h(D)>1$  as the density of primes represented by reduced forms of the discriminant $D$ is $1/2$. 

\begin{rmk}
 In the work of Matom$\ddot{\rm a}$ki  \cite[Lemma 6]{KMatomaki}, it was mentioned that the associated function, $\sigma(u)$ is the unique solution of the integral equation 
$$
u \sigma(u) = \int_{0}^{u} \sigma(t) \beta(u-t) dt
$$
with the initial condition $\sigma(u)= u^{\beta_{0}-1}$ for some  $u \in (0, x_{1}]$. For the step function $\beta$, we obtained that  $\sigma(1.6334)>0$.
Hence,  for $u = u_{0} = 1.6334=  1/0.61222$, we have
\begin{equation}\label{PBound}
\begin{split}
 \sum_{n \le Y^{u_{0}} \atop \gcd(n,N) =1 }   h_{Y}(n) r^{*}_{D}(n) &= (\sigma(u_{0}) + o_{\beta, U}(1)) \frac{P(1) L(1, \chi_{D})^{\beta_{0}}}{\Gamma(\beta_{0})} (\log Y)^{\beta_{0}-1} Y^{u_{0}} >0.
\end{split}
\end{equation}
\end{rmk}
Following the work of Matom$\ddot{\rm a}$ki \cite{KMatomaki}, these computations for the integrals and $u_{0}$ have been done using the mathematical software Mathematica \cite{Mathematica}. Here, we provide the source code for these computation.

\begin{lstlisting}
K1 = 25; 
m = 5;

xk[k_] := If[k == 0, 0, 1/(K1 - k + 1)];

chi[k_] := (3 - 4*(Sin[Pi/(K1 - k + 2)])^2);

dsol[K1_] := 
  NDSolve[{u*d'[u] + (1 - chi[0]) d[u] + 
      Sum[d[u - xk[k]] (chi[k - 1] - chi[k]), {k, 1, K1}] == 0, 
    d[u /; u <= xk[1]] == If[u <= 0, 0, u^(chi[0] - 1)]}, 
   d, {u, 0, 3}];

f[u_] := Evaluate[d[u]] /. First[Evaluate[dsol[K1]]];
\end{lstlisting}

\noindent
\subsection{An estimate for $n_{sym^{2}f, D}$:}   In this part, we evaluate sub-convexity estimate for $n_{sym^{2}f, D}$.
 
\bigskip
\noindent 
\textbf{Proof of \thmref{ExtMatKLSW}:}
By definition,  $n_{sym^{2}f, D}$ be the greatest integer among all $n$ such that  $\lambda_{sym^{2}f}(n) \ge 0$ for every square-free positive integer $n <  n_{sym^{2}f, D},$ $\gcd(n,N) =1$ and  $n= \mathcal{Q}(\underline{x})$  for some $ \mathcal{Q} \in \mathcal{S}_{D}$ and  $ \underline{x} \in \mathbb{Z}^{2}$.  Let us write  $n_{sym^{2}f, D} = Y = X^{0.61222}.$ 
 A key  idea to estimate $Y$ is to compare a  positive lower bound and upper bound for the sum
\begin{equation*}
\begin{split}
S(sym^{2}f, D; X ) 
& = \sideset{}{^{\flat }} \sum_{n \le X \atop \gcd(n,N) =1  } \lambda_{sym^{2}f}(n) r_{Q}(n).\\
\end{split}
\end{equation*}
From \thmref{PropUpper} and \propref{PropLower}, we have (for $X = Y^{1/0.61222}$)
\begin{equation*}
\begin{split}
 0< X\log^{2} X L(1, \chi_{D})^{3} \ll \sideset{}{^{\flat }} \sum_{n \le X \atop \gcd(n,N) =1 } h_{Y}(n) r^{*}_{D}(n) \le S(sym^{2}f, D; X )   \ll {(N^{2}k^{2}  |D|^{\frac{3}{2}} )}^{\frac{1}{3}+\epsilon} X^{\frac{2}{3}+\epsilon}.
\end{split}
\end{equation*} 
Thus, we have
\begin{equation*}
\begin{split}
 X^{1+\epsilon}   \ll {(N^{2}k^{2}  |D|^{\frac{3}{2}} )}^{\frac{1}{3}+\epsilon} L(1, \chi_{D})^{-3}  X^{\frac{2}{3}+\epsilon} \qquad \Longrightarrow  \qquad X \ll {(N^{2}k^{2}  |D|^{\frac{3}{2}})}^{1+\epsilon} L(1, \chi_{D})^{-9} .
 \end{split}
\end{equation*} 
Now, we put $X = Y^{1/0.61222}$ to get 
$$
Y \ll {(N^{2}k^{2}  |D|^{\frac{3}{2}})}^{0.61222 + \epsilon} L(1, \chi_{D})^{-5.51}.
$$
This shows that  $n_{{sym^{2}f, D}} \ll {(N^{2}k^{2}  |D|^{\frac{3}{2}} )}^{0.61222 + \epsilon}L(1, \chi_{D})^{-5.51}. $ 
This completes the proof.

\bigskip
\textbf{Acknowledgement :}  
The second author would like to thank Stat-Math Unit, ISI Delhi for giving opportunity to explore research activity. 


\smallskip
\textbf{Declaration:} There is no conflict of interest between authors, and to anyone.


\end{document}